\documentclass[10pt,a4paper]{article}

\usepackage{mathptmx}
\usepackage{amsmath}
\usepackage{amsfonts}  
\usepackage{amssymb}
\usepackage{amsthm}
\usepackage{setspace}
\usepackage{enumitem}
\usepackage{times}
\usepackage{relsize}
\usepackage{mathrsfs}
\usepackage{scalerel}

\DeclareMathAlphabet{\mathcal}{OMS}{cmsy}{m}{n}

\DeclareSymbolFont{dgmetaletters}{OML}{DGMetaSerifScience-TLF}{m}{it}
\DeclareMathSymbol{v}{\mathalpha}{dgmetaletters}{`v}
\DeclareMathSymbol{w}{\mathalpha}{dgmetaletters}{`w}
\DeclareMathSymbol{u}{\mathalpha}{dgmetaletters}{`u}
\DeclareMathSymbol{x}{\mathalpha}{dgmetaletters}{`x}

\usepackage[a4paper, left=4.2cm, right=4.2cm, top=4.4cm, bottom=4.4cm]{geometry}
\usepackage{setspace}

\allowdisplaybreaks

\makeatletter
\def\blfootnote{\xdef\@thefnmark{}\@footnotetext}
\makeatother

\newcommand{\dom}{\text{dom}}
\newcommand{\doititle}{\textit}

\newcommand{\newalpha}{\scaleobj{0.9}{\alpha}}
\newcommand{\newbeta}{\scaleobj{0.9}{\beta}}
\newcommand{\newmu}{\scaleobj{0.9}{\mu}}
\newcommand{\newnu}{\scaleobj{0.9}{\nu}}
\newcommand{\newtheta}{\scaleobj{0.9}{\theta}}

\newcommand{\newsigma}{\scaleobj{0.9}{\sigma}}
\newcommand{\newomega}{\scaleobj{0.9}{\omega}}

\newcommand{\newgamma}{\scaleobj{0.9}{\gamma}}
\newcommand{\newtau}{\scaleobj{0.9}{\tau}}
\newcommand{\neweta}{\scaleobj{0.9}{\eta}}

\newcommand{\infinity}{\mathrel{\raisebox{-.3ex}{$\mathlarger{\infty}$}}}

\newcommand{\textx}{x}
\newcommand{\mathx}{x}
\newcommand{\textv}{v}
\newcommand{\mathv}{v}
\newcommand{\textu}{u}
\newcommand{\mathu}{u}
\newcommand{\textw}{w}
\newcommand{\mathw}{w}

\numberwithin{equation}{section}

\newtheoremstyle{noenddot}
{\topsep}
{\topsep}%
{\itshape}
{-1pt}%
{\bfseries}
{}%
{ } 
{\thmname{#1}\thmnumber{ #2}\thmnote{ \normalfont(#3)}}%

\theoremstyle{noenddot}

\newtheorem{theorem}{Theorem}[section]

\newtheorem{lemma}[theorem]{Lemma}

\theoremstyle{definition}

\begin{document}
	
	\renewcommand{\abstractname}{\vspace{-\baselineskip}}
	
	\begin{center}
		\begin{spacing}{1.5}	
			\textbf{\huge{Fiber entropy and algorithmic complexity of random orbits}}
		\end{spacing}
		$~$
		
		\large{\textsc{Elias Zimmermann}}\blfootnote{The author is supported by GIF grant I-1485-304.6/2019.}
	\end{center}
	
	\begin{center}
		Mathematical Institute, University of Leipzig
		
		Augustusplatz 10, 04109 Leipzig
		
		elias.zimmermann@math.uni-leipzig.de
	\end{center}
	
	\begin{abstract} 
		\noindent ABSTRACT. Let $\Theta$ be a finite alphabet. We consider a bundle of measure preserving transformations $(T_{\newtheta})_{\newtheta \in \Theta}$ acting on a probability space $(X,\newmu)$, which are chosen randomly according to an ergodic stochastic process $(\Xi,\newnu,\newsigma)$ with state space $\Theta$. This describes a paradigmatic case of a random dynamical system (RDS). Considering a finite partition $\mathcal{P}$ of $X$ we show that the conditional algorithmic complexity of a random orbit $x, T_{\newalpha_{0}}(x),T_{\newalpha_{1}}\circ T_{\newalpha_{0}}(x),...$ in $X$ along a sequence $\newalpha = \newalpha_{0}\newalpha_{1}\newalpha_{2}...$ in $\Xi$ equals almost surely the fiber entropy of the RDS with respect to $\mathcal{P}$, whenever the latter is ergodic. This extends a classical result of A. A. Brudno connecting algorithmic complexity and entropy in deterministic dynamical systems.
	\end{abstract}
	
	\section{Introduction}
	
	In the theory of algorithmic complexity the randomness of words over finite alphabets is quantified in terms of their Kolmogorov complexity, which measures the length of a minimal description with respect to a universal Turing machine. Based on this notion the upper (lower) algorithmic complexity of an infinite sequence over a finite alphabet may be defined as the limes superior (inferior) of the normalized  Kolmogorov complexities of its initial segments. Given a stationary stochastic process with finite state space a classical theorem going back to Brudno states that upper and lower algorithmic complexity of a typical realization coincide and equal the entropy rate of the process, whenever the latter is ergodic.
	
	Brudno's theorem also admits a version for general dynamical systems, which can be obtained as follows. Let $(X,\newmu,T)$ be a measure preserving dynamical system (MDS) and $\mathcal{P}$ be a finite partition of $X$. Then to $\newmu$-almost every state $\textx \in X$ we may associate a unique \textit{$\mathcal{P}$-name}, i.\@ e. a sequence $\newomega \in \mathcal{P}^{\mathbb{N}}$ satisfying $T^{n}(\textx) \in \newomega_{n}$ for all $n \in \mathbb{N}$. The \textit{upper} and \textit{lower orbit complexity} $\overline{\mathcal{C}}_{\mathcal{P}}(\textx)$ and $\underline{\mathcal{C}}_{\mathcal{P}}(\textx)$ of $\textx$ with respect to $\mathcal{P}$ may then be defined as the upper and lower algorithmic complexity of $\newomega$. Denoting by $h_{\newmu}(\mathcal{P})$ the entropy rate of the transformation $T$ with respect to $\mathcal{P}$ Brudno's theorem reads as follows.
	
	\begin{theorem} \label{Brudno}
		\textnormal{\textbf{(Brudno, '82)}} Let $(X,\newmu,T)$ be an ergodic MDS and $\mathcal{P}$ be a finite partition of $X$. Then for $\newmu$-almost all $\mathx \in X$ we have
		\[\overline{\mathcal{C}}_{\mathcal{P}}(\mathx) = \underline{\mathcal{C}}_{\mathcal{P}}(\mathx) = h_{\newmu}(\mathcal{P}).\]
	\end{theorem}
	
	In \cite{Bru82} Brudno explicitly states the above identity only for the upper complexity. However, his proof can be modified in a well known way to obtain the same statement for the lower complexity. Moreover, Brudno proved his theorem only for symbolic shifts. Yet a standard argument involving the construction of a symbolic factor leads to the above version for general dynamical systems.
	
	Brudno's theorem is a remarkable statement for two reasons: On the one hand it connects two concepts of randomness (entropy and algorithmic complexity), which are defined in rather different ways, in a very close manner. On the other hand it provides a rigorous justification for the interpretation of entropy as a measure of orbit complexity.
	
	Generalizations of Brudno's theorem have been proposed for quantum statistics, see \cite{Ben07} and \cite{BKM+06}, and infinite measure spaces, see \cite{Zwe06}. More recently, several extensions of growing generality for actions of amenable groups were obtained. While a version for $\mathbb{Z}^{d}$-subshifts is proposed in \cite{FT17}, an extension to amenable groups admitting certain computable F\o lner monotilings is provided in \cite{Mor21}. Finally, a generalization for arbitrary computable amenable groups is proved in \cite{Alp18}. For connections of entropy and algorithmic complexity in topological dynamical systems see \cite{Alp18}, \cite{Bru82}, \cite{GHR10}, \cite{Mor15}, \cite{Sim15}, \cite{Whi91} and \cite{Whi93}.
	
	In this paper we propose a generalization of Brudno's result to the context of random dyna\-mics. The random dynamical systems we consider are given by a bundle of measure preserving transformations $(T_{\newtheta})_{\newtheta \in \Theta}$ over some finite alphabet $\Theta$, which act on a common probability space $(X,\newmu)$ and are chosen randomly according to an ergodic stochastic process $(\Xi,\newnu,\newsigma)$ with state space $\Theta$. This setting gives rise to a step skew product $T$ on the product space $\Xi \otimes X$, which preserves the product measure $\newnu \otimes \newmu$. The corresponding MDS $(\Xi \otimes X, \newnu \otimes \newmu, T)$ defines an instance of a \textit{bundle random dynamical system (bundle RDS)}, to which we shall refer as a \textit{finite bundle RDS}. A detailed introduction is given in Section \ref{FibEntr}. As we shall see, the above framework allows the formulation of a quite natural generalization of Brudno's theorem (which seems elusive in the setting of general random dynamical systems).
	
	To obtain such a generalization we shall use a conditional version of Kolmogorov complexity, which can be defined in terms of so called oracle machines. Based on this notion we will introduce the \textit{upper} and \textit{lower conditional algorithmic complexity} of an infinite sequence $\newomega$ relative to some oracle $\newalpha$, which will be given by another infinite sequence, as the limes superior (inferior) of the normalized conditional complexities of $\newomega$'s initial segments relative to $\newalpha$. For more details see Section \ref{Prel}.
	
	Considering a finite bundle RDS $(\Xi \otimes X, \newnu \otimes \newmu, T)$ together with a finite partition $\mathcal{P}$ of $X$ and fixing some sequence $\newalpha \in \Xi$ we may associate to $\newmu$-almost every state $\textx \in X$ a unique sequence $\newomega$ over $\mathcal{P}$ satisfying $T_{\newalpha_{n-1}} \circ ... \circ T_{\newalpha_{0}} (\textx) \in \newomega_{n}$ for all $n \in \mathbb{N}$, which we shall call the \textit{$\mathcal{P}$-$\newalpha$-name} of $\textx$. This allows us to introduce the \textit{upper} and \textit{lower conditional orbit complexity} $\overline{\mathcal{C}}_{\mathcal{P}}(\textx|\newalpha)$ and $\underline{\mathcal{C}}_{\mathcal{P}}(\textx|\newalpha)$
	of $\textx$ along $\newalpha$ as the upper and lower conditional algorithmic complexity of $\newomega$ relative to $\newalpha$.
	
	The entropy of an RDS is usually identified with the so called fiber entropy of the corresponding skew product, which will be defined rigorously in Section \ref{FibEntr}. Denoting by $h_{\newnu,\newmu}(\mathcal{P})$ the fiber entropy of $T$ with respect to the partition $\mathcal{P}$ we may formulate our main result as follows.
	
	\begin{theorem} \label{main}
		Let $(\Xi \times X,\newnu \otimes \newmu, T)$ be an ergodic finite bundle RDS and $\mathcal{P}$ be a finite partition of $X$. Then for $\newnu \otimes \newmu$-almost all $(\newalpha,x)$ we have
		\[\overline{\mathcal{C}}_{\mathcal{P}}(x|\newalpha) = \underline{\mathcal{C}}_{\mathcal{P}}(x|\newalpha) = h_{\newnu,\newmu}(\mathcal{P}).\]
	\end{theorem}
	Furthermore, denoting by $\mathcal{N}$ the natural partition of $\Xi$ consisting of the cylinder sets, which are only specified in the first symbol, we shall obtain the following decomposition formula for orbit complexities, which resembles the classical Abramov-Rokhlin decomposition of the fiber entropy of a skew product. \newpage
	
	\begin{theorem} \label{AR-comp}
		Let $(\Xi \times X, \newnu \otimes \newmu, T)$ be an ergodic finite bundle RDS and $\mathcal{P}$ be a finite partition of $X$. Then for $\newnu\otimes\newmu$-almost all $(\newalpha,\mathx)$ we have $\overline{\mathcal{C}}_{\mathcal{N}\times\mathcal{P}}(\newalpha,\mathx) = \underline{\mathcal{C}}_{\mathcal{N}\times\mathcal{P}}(\newalpha,\mathx)$, $\overline{\mathcal{C}}_{\mathcal{P}}(\mathx|\newalpha) = \underline{\mathcal{C}}_{\mathcal{P}}(\mathx|\newalpha)$ and $\overline{\mathcal{C}}(\newalpha) = \underline{\mathcal{C}}(\newalpha)$ and, denoting by ${\mathcal{C}}_{\mathcal{N} \times \mathcal{P}}(\newalpha,\mathx)$, ${\mathcal{C}}_{\mathcal{P}}(\mathx|\newalpha)$ and $\mathcal{C}(\newalpha)$ the respective common values, obtain the relation
		\[\mathcal{C}_{\mathcal{P}}(\mathx|\newalpha) = \mathcal{C}_{\mathcal{N} \times \mathcal{P}}(\newalpha,\mathx) - \mathcal{C}(\newalpha).\]
	\end{theorem}
	
	The paper is organized as follows: Section \ref{Prel} contains the necessary preliminaries from the theory of algorithmic complexity and the entropy theory of deterministic systems. In Section \ref{FibEntr} we introduce the relevant aspects of entropy theory for finite bundle RDS' and apply them to certain examples arising from actions of finitely generated groups. Finally, Section \ref{Proofs} is devoted to the proof of our main results.
	
	\section{Preliminaries} \label{Prel}
	
	We set $\mathbb{N} = \{0,1,2,...\}$. Throughout this paper an alphabet is always a finite, non empty set. Given an alphabet $\Lambda$ we shall write $\Lambda^{*}$ for the set of (finite) words over $\Lambda$. The empty word is denoted by $\mathbf{e}$. We say that a word $\textv \in \Lambda^{*}$ is a \textit{prefix} of a word $\textw \in \Lambda^{*}$ in case $|\textv| \le |\textw|$ and $\textv_{i} = \textw_{i}$ for $i \in \{0,...,|\textv|{-}1\}$. By convention the empty word is a prefix of every word. A set $B \subseteq \Lambda^{*}$ is called \textit{prefix free} in case $\textv = \textw$ for any $\textv,\textw \in B$ such that $\textv$ is a prefix of $\textw$.
	
	\subsection{Kolmogorov complexity}
	
	Kolmogorov complexity is usually defined in terms of Turing machines. A \textit{Turing machine} is a mathematical model of a computer. Informally, it
	consists of three different tapes built of cells, which are either blank or inscribed with a symbol, and equipped with devices for reading and writing, the so called \textit{heads}. In particular, there is a right-infinite tape together with a reading head, the so called \textit{input tape}, a bi-infinite tape together with a head for reading and writing, the \textit{work tape}, and a right-infinite tape together with a writing head, the so called \textit{output tape}. At the beginning of a computation the input word is written on the input tape,
	while the machine is in the start state. The computation proceeds as follows. In every step
	the heads on the input and work tapes read the cell at their actual position. Depending on the scanned inscriptions and the actual state of the machine the head on the work tape may overwrite the actual cell and move to the left
	or to the right, while the head on the output tape stands still or prints a symbol at its actual position and moves one cell to the right. Finally the head on the input tape moves one cell to the right as the machine changes into a new state. All this happens in accordance
	with a fixed program. If a so called stop state is reached, the machine halts and outputs
	the word written on the output tape. The machine may also never reach a stop
	state, in which case it computes infinitely long.
	
	There are various equivalent formalizations of Turing machines. However, following a common practice, we shall confine ourselves to the above informal description. For further information and references on Turing computability we refer the reader to \cite{DH10}. In the following we shall always consider Turing machines $M$ with input alphabet $\{0,1\}$ and some fixed output alphabet $\Lambda$. Given such a machine we shall write $\dom(M)$ for the set of words in $\{0,1\}^{*}$, on which $M$ terminates. For an input $\textu \in \dom(M)$ we shall denote the output of $M$ corresponding to $\textu$ by $M(\textu)$. We will call a machine $M$ \textit{prefix free} in case $\dom(M)$ is a prefix free set. By a \textit{universal (prefix free)} Turing machine we mean a (prefix free) Turing machine $U$, which can simulate every other (prefix free) Turing machine when given its programme in a suitable form as additional input. More precisely we require that for every (prefix free) machine $M$ there is a word $\textw \in \{0,1\}^{*}$ such that $M$ terminates on an input $\textu$ if and only if $U$ terminates on the input $\textw \textu$ and the the outputs $M(\textu)$ and $U(\textw \textu)$ coincide in this case. The existence of a universal (prefix free) Turing machine is a well known fact of computability theory, see e.\@ g. \cite[Ch. 2 and 3]{DH10}.
	
	Kolmogorov complexity as well as conditional Kolmogorov complexity, which is defined later, are usually only considered for words (and oracles) over $\{0,1\}$. However, it will be convenient for us to define both for words (and oracles) over an arbitrary alphabet. The properties proposed in this and the next section are obtained by a straightforward adaption of the standard proofs, as e.\@ g. given in \cite{DH10}, to this slightly extended setting.
	
	For a given Turing machine $M$ we define the \textit{complexity} of a word $\textv \in \Lambda^{*}$ with respect to $M$ as the number
	\[C_{M}(\mathv) := \min\big\{|\mathu|\colon \mathu \in \dom(M)\text{ and }M(\mathu) = \mathv\big\},\]
	where $\min\emptyset = \mathrel{\raisebox{-.1ex}{$\mathlarger{\mathlarger{\infty}}$}}$ by convention. In other words $C_{M}(\textv)$ measures the length of a minimal description of $\textv$ with respect to $M$. To obtain a measure of complexity, which takes every possible choice of $M$ into account, we fix a universal Turing machine $U$ and define the \textit{(Kolmogorov) complexity} $C(\textv)$ of $\textv$ as the complexity $C_{U}(\textv)$ of $\textv$ with respect to $U$. By the defining property of a universal machine we obtain then
	\[C(\textv) \le C_{M}(\textv) + O(1)\]
	for any Turing machine $M$, where the constant $O(1)$ depends only on $M$. In parti\-cular, choosing another universal machine in the definition will not change the values of $C$ up to an additive constant. In this sense $C$ is independent under the choice of $U$, which is the commonly accepted justification for the above definition.
	
	An often considered variant of $C$ is obtained by using a universal prefix free Turing machine $P$ instead of $U$. The corresponding complexity $C_{P}(\textv)$, which is denoted by $K(\textv)$, is called the \textit{prefix free complexity} of $\textv$. Similarly as above one has
	\[K(\mathv) \le C_{M}(\mathv) + O(1)\]
	for every prefix free Turing machine $M$ with some constant $O(1)$ depending only on $M$. In particular, $K$ is independent under the choice of the machine $P$ up to an additive constant. The reason for considering $K$ instead of $C$ is that the former has some technical advantages in comparison to the latter. Most importantly it satisfies the subadditivity property
	\[K(\mathv \mathw) \le K(\mathv) + K(\mathw) + O(1)\]
	for all $\textv,\textw \in \Lambda^{*}$ and some constant $O(1)$, which is useful in many situations. However, the difference between $K$ and $C$ is bounded by
	\begin{align} \label{Equ1}
		|C(\mathv) - K(\mathv)| \le 2\log|\mathv| + O(1)
	\end{align}
	for all $\textv \in \Lambda^{*}$ and some constant $O(1)$. Furthermore, considering a Turing machine, which maps the $n$-th word in $\{0,1\}^{*}$ to the $n$-th word in $\Lambda^{*}$, one obtains
	\begin{align} \label{Equ2}
		C(\mathv) \le |\mathv|\log|\Lambda| + O(1)
	\end{align}
	for every $\textv \in \Lambda^{*}$ and some constant $O(1)$.
	
	The upper bound on the difference of $C$ and $K$ given in (\ref{Equ1}) makes it irrelevant, which complexity measure is used when passing to complexity rates of infinite sequences. More precisely, denoting by $C_{n}(\newomega)$ and $K_{n}(\newomega)$ the complexities $C(\newomega_{0}...\newomega_{n-1})$ and $K(\newomega_{0}...\newomega_{n-1})$ one obtains
	\[\limsup_{n \to \infinity} \frac{\mathrel{\raisebox{-.4ex}{$1$}}}{\mathrel{\raisebox{+.1ex}{$n$}}}C_{n}(\newomega) = \limsup_{n \to \infinity} \frac{\mathrel{\raisebox{-.4ex}{$1$}}}{\mathrel{\raisebox{+.1ex}{$n$}}}K_{n}(\newomega)\]
	for every infinite sequence $\newomega$ over $\Lambda$. We shall call the above value the \textit{upper algorithmic complexity} of $\newomega$ and denote it by $\overline{\mathcal{C}}(\newomega)$. Of course the analogous statement holds for the limes inferior, which we shall denote by $\underline{\mathcal{C}}(\newomega)$ and call the \textit{lower algorithmic complexity} of $\newomega$. By (\ref{Equ2}) it is clear that $\overline{\mathcal{C}}(\newomega)$ and $\underline{\mathcal{C}}(\newomega)$ take values in $[0,\log|\Lambda|]$. Recall that the difference resulting from underlying different universal Turing machines in the definition of $K$ and $C$ is bounded by a constant, so the arbitrariness in the definition of $K$ and $C$ discussed earlier does not matter at all for the definition of $\overline{\mathcal{C}}(\newomega)$ and $\underline{\mathcal{C}}(\newomega)$.
	
	\subsection{Conditional complexity} \label{Prel2}
	
	For the generalization of Brudno's theorem proposed in Theorem  \ref{main} we need an extension of the concepts discussed so far. More specifically we will require a conditional version of Kolmogorov complexity. Such a version can be obtained in terms of oracle machines. Informally an \textit{oracle machine} is a Turing machine with an additional one-sided tape, the so called \textit{oracle tape}, which is inscribed with an infinite sequence, the so called \textit{oracle}, and equipped with a reading head. In every step of the computation this head may read the inscription of the oracle tape at its actual position and move one cell to the right. Apart from that the operation principle is analogous to that of an ordinary Turing machine except for the fact that the actions on the work and output tape are executed not only in dependence of the symbol read on the input and work tape and the current state of the machine, but also on the symbol read on the oracle tape.
	
	As before we shall confine ourselves to the above informal description of oracle machines and refer the reader to \cite{DH10} for more information and literature. Again we shall only consider oracle machines $M$ with input alphabet $\{0,1\}$ and certain fixed output and oracle alphabets $\Lambda$ and $\Theta$. Given such an oracle machine we shall write $\dom(M,\newalpha)$ for the set of words, on which $M$ terminates when provided with $\newalpha$ as oracle, and $M^{\newalpha}(\textu)$ for the output of $M$ corresponding to an input $\textu$ out of $\dom(M,\newalpha)$. We shall call $M$ \textit{prefix free} if for every oracle $\newalpha$ the set $\dom(M,\newalpha)$ is a prefix free set. By a \textit{universal (prefix free)} oracle machine we mean a (prefix free) oracle machine $U$, which can simulate every other (prefix free) oracle machine $M$ when provided with its programme in a suitable form. More precisely we require that for every such machine $M$ there is a word $\textw \in \{0,1\}^{*}$ such that for all oracles $\newalpha$ we have $\textu \in \dom(M,\newalpha)$ if and only $\textw \textu \in \dom(U,\newalpha)$ and the outputs $M^{\newalpha}(\textu)$ and $U^{\newalpha}(\textw \textu)$ coincide in this case. The existence of universal (prefix free) oracle machines is also well known, see \cite[Ch. 2 and 3]{DH10}.
	
	Based on oracle machines a conditional variant of Kolmogorov complexity can be defined as follows. Given an oracle machine $M$ provided with an oracle $\newalpha$ one defines the \textit{conditional complexity} of a word $\textv \in \Lambda^{*}$ with respect to $M$ relative to $\newalpha$ as the number
	\[C_{M}(\mathv|\newalpha) := \min\big\{|\mathu|\colon \mathu \in \dom(M,\newalpha)\text{ and }M^{\newalpha}(\mathu) = \mathv\big\}.\]
	As in the ordinary case the \textit{conditional complexity} $C(\textv|\newalpha)$ of a word $\textv$ relative to $\newalpha$ is then defined as the conditional complexity $C_{U}(\textv|\newalpha)$ with respect to some fixed universal oracle machine $U$. Using a prefix free universal  oracle machine $P$ instead of $U$ we obtain the \textit{prefix free conditional complexity} $K(\textv|\newalpha)$ of $\textv$ relative to $\newalpha$. As in the ordinary case the definitions are invariant under the choice of the underlying universal oracle machines up to an additive constant. The latter is a consequence of the fact that we have
	\begin{equation} \label{Equ3}
		C(\textv|\newalpha) \le C_{M}(\textv|\newalpha) + O(1)
	\end{equation}
	for every oracle machine $M$ with a constant $O(1)$ depending only on $M$. The analogous statement is true for the prefix free conditional complexity, where $M$ ranges over all prefix free oracle machines. Furthermore, in the prefix free case we have
	\begin{equation} \label{Equ4}
		K(\mathv\mathw|\newalpha) \le K(\mathv|\newalpha) + K(\mathw|\newalpha) + O(1)
	\end{equation}
	for all $\textv,\textw \in \Lambda^{*}$, while the difference of $C(\textv|\newalpha)$ and $K(\textv|\newalpha)$ is again bounded from above by $2\log|\textv| + O(1)$ for some constant $O(1)$. In combination with (\ref{Equ4}) this implies
	\begin{equation} \label{Equ5}
		C(\mathv\mathw|\newalpha) \le C(\mathv|\newalpha) + C(\mathw|\newalpha) + 2\log|\mathv||\mathw| + O(1)
	\end{equation}
	for all $\textv,\textw \in \Lambda^{*}$ and some constant $O(1)$.
	
	As a further consequence it makes again no difference, which measure is used when considering conditional complexity rates of infinite sequences. More precisely, denoting by $C_{n}(\newomega|\newalpha)$ and $K_{n}(\newomega|\newalpha)$ the complexities $C(\newomega_{0}...\newomega_{n-1}|\newalpha)$ and $K(\newomega_{0}...\newomega_{n-1}|\newalpha)$ one obtains
	\[\limsup_{n \to \infinity} \frac{\mathrel{\raisebox{-.4ex}{$1$}}}{\mathrel{\raisebox{+.1ex}{$n$}}}C_{n}(\newomega|\newalpha) = \limsup_{n \to \infinity} \frac{\mathrel{\raisebox{-.4ex}{$1$}}}{\mathrel{\raisebox{+.1ex}{$n$}}}K_{n}(\newomega|\newalpha)\]
	for every infinite sequence $\newomega$ and every oracle $\newalpha$. We shall denote the above value by $\overline{\mathcal{C}}(\newomega|\newalpha)$ and call it the \textit{upper conditional algorithmic complexity} of $\newomega$ relative to $\newalpha$. The analogous statement is true for the limes inferior, which we shall denote by $\underline{\mathcal{C}}(\newomega|\newalpha)$ and call the \textit{lower conditional algorithmic complexity} of $\newomega$ relative to $\newalpha$. Note that as before the values are independent of the choice of the underlying universal oracle machines.
	
	It is not difficult to verify that the conditional complexity $C(\textv|\newalpha)$ of $\textv$ relative to $\newalpha$ is bounded from above by the plain complexity $C(\textv)$ of $\textv$ modulo some constant. Consequently the conditional upper (lower)  algorithmic complexity of a sequence is always bounded by its (unconditional) upper (lower) algorithmic complexity. This implies that $\overline{\mathcal{C}}(\newomega|\newalpha)$ and $\underline{\mathcal{C}}(\newomega|\newalpha)$ take again values in $[0,\log|\Lambda|]$. Furthermore one has $\overline{\mathcal{C}}(\newomega) = \overline{\mathcal{C}}(\newomega|\newtheta^{\infty})$ and $\underline{\mathcal{C}}(\newomega) = \underline{\mathcal{C}}(\newomega|\newtheta^{\infty})$ for every $\newtheta \in \Theta$, where $\newtheta^{\infty}$ denotes the infinite concatenation of $\newtheta$. Therefore algorithmic complexities arise from conditional algorithmic complexities as a special case. Finally it is not difficult to check that the values of the upper conditional algorithmic complexity are independent of finite initial segments of $\newomega$ and $\newalpha$. More precisely we have
	\begin{equation} \label{Equ6}
		\overline{\mathcal{C}}(\mathv\newomega|\mathu\newalpha) = \overline{\mathcal{C}}(\newomega|\newalpha)
	\end{equation}
	for all words $\textv \in \Lambda^{*}$ and $\textu \in \Theta^{*}$. The same is true for the lower conditional algorithmic complexity. This invariance property will become important below.
	
	\subsection{Measure theoretic entropy}
	
	Let $(X,\newmu,T)$ be a measure preserving dynamical system (MDS) consisting of a standard probability space $(X,\scaleobj{0.9}{\mathcal{A}},\newmu)$ together with a measure preserving transformation $T$ on $X$. In the following we will often omit the specification of the $\newsigma$-algebra $\scaleobj{0.9}{\mathcal{A}}$. By a \textit{finite partition} of $X$ we mean a finite collection of a.\@ s. disjoint measurable sets of positive measure such that their union covers the space $X$ up to a null set. Given a finite partition $\mathcal{P}$ of $X$ the \textit{Shannon entropy} of $\mathcal{P}$ is defined as
	\[H_{\newmu}(\mathcal{P}) := -\sum_{P \in \mathcal{P}}\newmu(P)\log\newmu(P).\]
	Note that this definition is independent of the base of the logarithm up to a factor. Throughout this paper we will always work with the binary logarithm, which is more convenient in the context of complexities than the natural logarithm.
	
	Given finite partitions ${\mathcal{P}}_{0},...,{\mathcal{P}}_{n-1}$ of $X$ we define their \textit{common refinement} as the finite partition
	\[\bigvee_{i=0}^{n-1}{\mathcal{P}}_{\mathrel{\raisebox{-.4ex}{$\mathsmaller{i}$}}} := \left\{\bigcap_{i=0}^{n-1}P_{\mathrel{\raisebox{-.4ex}{$\mathsmaller{i}$}}}\colon P_{\mathrel{\raisebox{-.4ex}{$\mathsmaller{i}$}}} \in {\mathcal{P}}_{\mathrel{\raisebox{-.4ex}{$\mathsmaller{i}$}}}\right\}.\]	
	Denoting by $T^{-i}(\mathcal{P})$ the finite partition $\{T^{-i}(P)\colon\! P \in\!\mathcal{P}\}$ and by
	\[H^{n}_{\newmu}(\mathcal{P}) := H_{\newmu}\left(\bigvee_{i=0}^{n-1}T^{-i}(\mathcal{P})\right)\]
	the Shannon entropies of the refinements of $\mathcal{P}$ under the iterations of $T$ a standard subadditivity argument yields that the limit
	\[h_{\newmu}(\mathcal{P}) := \lim_{n \to \infinity}\frac{\mathrel{\raisebox{-.4ex}{$1$}}}{\mathrel{\raisebox{+.1ex}{$n$}}}H_{\newmu}^{n}(\mathcal{P})\]
	exists and coincides with the infimum of the sequence. The value $h_{\newmu}(\mathcal{P})$ is called the \textit{entropy} of $T$ with respect to $\mathcal{P}$. The \textit{(Kolmogorov-Sinai) entropy} $h_{\newmu}$ of $T$ is then defined as the supremum of $h_{\newmu}(\mathcal{P})$ over all finite partitions $\mathcal{P}$ of $X$. By the Kolmogorov-Sinai theorem the supremum is attained if $\mathcal{P}$ is a so called \textit{generating partition} or \textit{generator} of $X$. Krieger's finite generator theorem guarantees the existence of such partitions for all ergodic, invertible transformations with finite entropy, see \cite[Thm. 4.2.3]{Dow11}. Accordingly Brudno's theorem implies that the entropy of those transformations quantifies the orbit complexity of a typical trajectory with respect to a generating partition.
	
	In case of an ergodic transformation the entropy $h_{\newmu}(\mathcal{P})$ with respect to a partition $\mathcal{P}$ can also be obtained as the limit of an individual information function. Given  a finite partition $\mathcal{Q}$ of $X$ the \textit{information function} $J_{\newmu}(\mathcal{Q})\colon X \to [0,\mathlarger{\mathlarger{\infty}}]$ with respect to $\mathcal{Q}$ is defined by
	\[J_{\newmu}(\mathcal{Q})(\textx) := -\log\newmu\big(\mathcal{Q}(\textx)\big),\] where $\mathcal{Q}(\textx)$ denotes the a.\@ s. unique element $Q \in \mathcal{Q}$ such that $\textx \in Q$. It is not difficult to see that $J_{\newmu}(\mathcal{Q})$ is a measurable function. Denoting by
	\[J_{\newmu}^{~\!\!\! n}(\mathcal{P}) := J_{\newmu}\left(\bigvee_{i=0}^{n-1}T^{-i}(\mathcal{P})\right)\]
	the information functions of the dynamical refinements of $\mathcal{P}$ under $T$ one obtains
	\[\lim_{n \to \infinity}\frac{\mathrel{\raisebox{-.4ex}{$1$}}}{\mathrel{\raisebox{+.1ex}{$n$}}}J_{\newmu}^{~\!\!\! n}(\mathcal{P})(\textx) = h_{\newmu}(\mathcal{P})\]
	for $\newmu$-almost all $\textx \in X$, whenever $T$ is ergodic, see \cite[Thm. 2.5]{Par69}. This is the statement of the famous Shannon-McMillan-Breiman (SMB) theorem.
	
	A class of transformations, which will play a particularly important role for us, are \textit{shifts}. For a given alphabet $\Theta$ we shall denote by $\Xi$ the space $\Theta^{\mathbb{N}}$ of infinite sequences over $\Theta$, which we may equip with the $\newsigma$-algebra generated by cylinder sets. The \textit{shift map} $\newsigma$ on $\Xi$, which maps a sequence $\newalpha_{0}\newalpha_{1}\newalpha_{2}...$ in $\Xi$ to the shifted sequence $\newalpha_{1}\newalpha_{2}\newalpha_{3}...$, defines then a measurable transformation. Given a $\newsigma$-invariant probability measure $\newnu$ on $\Xi$ we shall call the arising MDS $(\Xi,\newnu,\newsigma)$ a \textit{shift (system)} with state space $\Theta$. For a word $\textv \in \Theta^{*}$ we shall denote by
	\[[\textv] : = \{\textv_{0}\} \times ... \times \{\textv_{|v|-1}\} \times \Theta^{\mathbb{N}}\] 
	the cylinder set consisting of all sequences in $\Xi$ extending $\textv$. In the following we will often write $\newnu[\textv]$ to denote the value of the measure $\newnu$ on the set $[\textv]$. We shall call the partition $\mathcal{N}$ consisting of the sets $[\newtheta]$ for $\newtheta \in \Theta$ the \textit{natural partition} of $\Xi$. It is well known that $\mathcal{N}$ is a generator for $\newsigma$. In particular, by the Kolmogorov-Sinai theorem the entropy $h_{\newnu}$ of the shift coincides with the entropy rate $h_{\newnu}(\mathcal{N})$.
	
	\section{Entropy of a finite bundle RDS} \label{FibEntr}
	
	\subsection{Finite bundles}
	
	Let $\Theta$ be an alphabet and $(\Xi,\newnu,\newsigma)$ be an ergodic shift with state space $\Theta$. Furthermore, let $(X,\newmu)$ be a standard probability space and $(T_{\newtheta})_{\newtheta \in \Theta}$ be a bundle of $\newmu$-preserving transformations on $X$. Within this setting we may define a skew product $T$ on $\Xi \times X$ by
	\[T(\newalpha,\textx) := \big(\newsigma(\newalpha),T_{\newalpha_{0}}(\textx)\big)\]
	for $(\newalpha,\textx) \in \Xi \times X$. Skew products of this form are sometimes called \textit{step skew products}. It is not difficult to see that $T$ is measurable and preserves the product measure $\newnu \otimes \newmu$. The corresponding MDS $(\Xi \times X, \newnu \otimes \newmu,T)$ is an instance of a \textit{bundle random dynamical system (bundle RDS)} as defined in \cite{Arn98} and \cite{KL06}. We shall refer to it as a \textit{finite bundle RDS}. As noted earlier one may think of the underlying shift as a stationary ergodic stochastic process choosing one out of finitely many transformations in every step of time.
	
	Historically, finite bundle RDS' were among the first instances of random dynamical systems that have been stu\-died systematically, cf. \cite{UvN45}. They play an important role in the ergodic theory of finitely generated groups, where they have for instance been used to obtain weighted ergodic theorems as well as asymptotic equipartition properties, see \cite{Buf}, \cite{Buf00}, \cite{Gri99}, \cite{Kak51}, \cite{NP22+} and \cite{Ose65}.
	
	Given a bundle of measure preserving transformations $(T_{\newtheta})_{\newtheta \in \Theta}$ as above we shall call a measurable set $A \subseteq X$ \textit{invariant} if up to null sets we have $T^{-1}_{\newtheta}(A) = A$ for all $\newtheta \in \Theta$. The bundle is called \textit{ergodic} if every invariant set has trivial measure. It is not difficult to see that ergodicity of the skew product implies ergodicity of the bundle. In applications it is often useful to know whether the converse implication holds, i.\@ e. whether ergodicity of the bundle implies ergodicity of $T$. In the case of a Bernoulli measure the validity of this implication was shown by Kakutani in \cite{Kak51}. For Markov measures an important criterion for its validity is due to Bufetov, see \cite{Buf00}.
	
	Recall that a measure $\newnu$ on $\Xi$ is called a \textit{Markov measure} if there exists a row stochastic matrix $\Pi \in [0,1]^{\Theta \times \Theta}$ together with a probability vector $\pi \in [0,1]^{\Theta}$ such that we have
	\begin{equation} \label{Equ7}
		\newnu[\textv] = \pi(\textv_{0})~\Pi(\textv_{0},\textv_{1})~...~\Pi(\textv_{|v|-2},\textv_{|v|-1})
	\end{equation}
	for all words $\textv \in \Theta^{*}$. In turn, given a vector $\pi$ and a matrix $\Pi$ as above, there is a unique probability measure $\newnu$ on $\Xi$ satisfying (\ref{Equ7}) for all words $\textv \in \Theta^{*}$ as a consequence of Kolmogorov's extension theorem. As one can show, $\newnu$ is $\newsigma$-invariant if and only if we have $\pi^{T}\Pi = \pi^{T}$. Moreover, if $\pi \in (0,1)^{\Theta}$, then $\newnu$ is ergodic if and only if $\Pi$ is irreducible. \pagebreak
	
	Let $(\Xi \times X, \newnu \otimes \newmu, T)$ be a finite bundle RDS and $\newnu$ be a $\newsigma$-invariant ergodic Markov measure with transition matrix $\Pi$ and initial probability vector $\pi \in (0,1)^{\Theta}$. Assume that the matrix $\Pi^{T}\Pi$ is irreducible. Then a theorem of Bufetov, cf. \cite[Thm. 5]{Buf00}, states that the ergodicity of the skew product is equivalent to the ergodicity of the bundle. A stronger condition requires that for all $\newtheta,\newtheta' \in \Theta$ there is some $\delta \in \Theta$ such that $\Pi(\delta,\newtheta) > 0$ and $\Pi(\delta,\newtheta') > 0$, which is of course satisfied if $\newnu$ is for instance a Bernoulli measure.
	
	An important instance of a finite bundle RDS is a \textit{random shift}. For alphabets $\Lambda$ and $\Theta$ let $\Omega$ denote the \textit{configuration space} $\Lambda^{(\Theta^{*})}$ consisting of families $\newbeta = (\newbeta_{u})_{u \in \Theta^{*}}$ with index set $\Theta^{*}$ and entries in $\Lambda$ and equip it with the $\newsigma$-algebra generated by cylinder sets. For $\newtheta \in \Theta$ let $S_{\newtheta}$ denote the \textit{$\newtheta$-shift} on $\Omega$ given by
	\[S_{\newtheta}(\newbeta) := (\newbeta_{\newtheta u})_{u \in \Theta^{*}}\]
	for $\newbeta \in \Omega$, which is obviously measurable. Finally, let $S$ denote the corresponding step skew product on $\Xi \times \Omega$, which takes the form
	\[S(\newalpha,\newbeta) = \big(\newsigma(\newalpha),S_{\newalpha_{0}}(\newbeta)\big)\]
	for $(\newalpha,\newbeta) \in \Xi \times \Omega$. Then the bundle RDS $(\Xi \times \Omega,\newnu \otimes \newmu,S)$ arising from any probability measure $\newmu$ on $\Omega$, which is invariant under all $S_{\newtheta}$, shall be called a \textit{random shift}.
	
	\subsection{Fiber entropy}
	
	The entropy of a finite bundle RDS $(\Xi \times X, \newnu \otimes \newmu, T)$ can be defined as follows. Let $\mathcal{P}$ be a finite partition of $X$. Writing $T_{\newalpha,n}$ for the transformation $T_{\newalpha_{n-1}} \circ ... \circ T_{\newalpha_{0}}$ and $T_{\newalpha,0}$ for $\text{Id}$ we may denote by
	\[J_{\newmu}^{~\!\!\! n}(\mathcal{P})(\newalpha,\textx) := J_{\newmu}\left(\bigvee_{i=0}^{n-1}T_{\newalpha,i}^{-1}(\mathcal{P})\right)(\textx)\]
	the values of the information function and by
	\[H_{\newmu}^{n}(\mathcal{P})(\newalpha):= H_{\newmu}\left(\bigvee_{i=0}^{n-1}T_{\newalpha,i}^{-1}(\mathcal{P})\right)\]
	the entropy values of the dynamical refinements of $\mathcal{P}$ along a given sequence $\newalpha \in \Xi$. To obtain an entropy notion for $T$ we have to consider the averaged entropies
	\[H_{\newnu,\newmu}^{n}(\mathcal{P}) := \int H^{n}_{\newmu}(\mathcal{P})(\newalpha)~d\newnu(\newalpha).\]
	Denoting again by $\mathcal{N}$ the natural partition of $\Xi$ we observe that
	\[\left(\bigvee_{i=0}^{n-1}T^{-i}(\mathcal{N} \times \mathcal{P})\right)(\newalpha,\textx) = \left(\bigvee_{i=0}^{n-1}\newsigma^{-i}(\mathcal{N})\right)(\newalpha) \times  \left(\bigvee_{i=0}^{n-1}T_{\newalpha,i}^{-1}(\mathcal{P})\right)(\textx)\]
	for $\newnu \otimes \newmu$-almost all $(\newalpha,\textx) \in \Xi \times X$. Using this it is not difficult to verify that
	\begin{equation} \label{Equ8}
		J^{~\!\!\! n}_{\newmu}(\mathcal{P})(\newalpha,\textx) = J^{~\!\!\! n}_{\newnu\otimes\newmu}(\mathcal{N} \times \mathcal{P})(\newalpha,\textx) - J^{~\!\!\! n}_{\newnu}(\mathcal{N})(\newalpha)
	\end{equation}
	for $\newnu\otimes \newmu$-almost all $(\newalpha,\textx) \in \Xi \times X$. Via integration one obtains	therefore\pagebreak 
	\[H_{\newnu,\newmu}^{n}(\mathcal{P}) = H^{n}_{\newnu \otimes \newmu}(\mathcal{N} \times \mathcal{P}) - H_{\newnu}^{n}(\mathcal{N}),\]
	which implies that the limit
	\[h_{\newnu,\newmu}(\mathcal{P}) = \lim_{n \to \infinity}\frac{\mathrel{\raisebox{-.4ex}{$1$}}}{\mathrel{\raisebox{+.1ex}{$n$}}}H_{\newnu,\newmu}^{n}(\mathcal{P})\]
	exists and satisfies
	\begin{equation} \label{Equ9}
		h_{\newnu,\newmu}(\mathcal{P}) = h_{\newnu \otimes \newmu}(\mathcal{N} \times \mathcal{P}) - h_{\newnu}.
	\end{equation}
	The value $h_{\newnu,\newmu}(\mathcal{P})$ is called the \textit{fiber entropy (rate)} of $T$ with respect to $\mathcal{P}$. The \textit{fiber entropy} $h_{\newnu,\newmu}$ of $T$ is now defined as the supremum of $h_{\newnu,\newmu}(\mathcal{P})$ over all finite partitions $\mathcal{P}$ of $X$.  The decomposition in (\ref{Equ9}) may be seen as a partitionwise version of the well known Abramov-Rohlin formula
	\begin{align*}
		h_{\newnu,\newmu} = h_{\newnu \otimes \newmu} - h_{\newnu},
	\end{align*}
	which is true in a much more general setting and was first proved by Abramov and Rokhlin in \cite{AR62}. 
	
	Using (\ref{Equ8}) one can also deduce a random version of the SMB theorem. In fact, if $T$ is ergodic, the classical SMB theorem yields that the right hand side of (\ref{Equ8}) divided by $n$ converges $\newnu \otimes \newmu$-almost surely to the right hand side of (\ref{Equ9}) (note that $\newsigma$ is ergodic by assumption). Thus we obtain
	\[\lim_{n \to \infinity}\frac{\mathrel{\raisebox{-.4ex}{$1$}}}{\mathrel{\raisebox{+.1ex}{$n$}}}J_{\newmu}^{~\!\!\! n}(\mathcal{P})(\newalpha,\textx) = h_{\newnu,\newmu}(\mathcal{P})\]
	for $\newnu \otimes \newmu$-almost all $(\newalpha,\textx) \in \Xi \times X$ in this case.
	
	It should be noted that fiber entropy can be defined for much more general random dynamical systems as a certain relative (conditional) entropy of the corresponding skew product, see \cite{Bog90}, \cite{Bog93}, \cite{KL06} and \cite{Mor86}, containing the above notion of fiber entropy as a special case. A general version of the random SMB theorem is then obtained as an instance of the relative SMB theorem. The same way one can establish a random version of the Kolmogorov-Sinai theorem, which implies that $h_{\newnu,\newmu}$ coincides with $h_{\newnu,\newmu}(\mathcal{P})$ whenever $\mathcal{P}$ is a \textit{random generator}. The latter means that the smallest $\newsigma$-algebra containing the partitions $T_{\newalpha,n}^{-1}(\mathcal{P})$ for all $n \in \mathbb{N}$ equals $\newnu$-almost surely the whole $\newsigma$-algebra on $X$ (up to null sets). Accordingly, by the random version of Brudno's Theorem proposed in Theorem \ref{main}, the fiber entropy of a finite bundle RDS admitting a random generator equals almost surely the conditional complexity of a random orbit with respect to this generator.
	
	\subsection{Examples}
	
	In the remaining part of this section we shall discuss two examples of ergodic finite bundle RDS' arising from actions of finitely generated groups and compute their fiber entropy. These examples will exhaust the range of values fiber entropy can take. Let $G$ be a countable group with neutral element $e$ and $(X,\newmu)$ be a standard probability space. By a \textit{measure preserving action} of $G$ on $X$ we mean a family $\{T_{g}\colon g \in G\}$ of measure preserving invertible transformations $T_{g}$ on $X$ such that $T_{e} = \text{Id}$ and $T_{g} \circ T_{h} = T_{gh}$ for $g,h \in G$. A measurable set $A \subseteq X$ is called \textit{$G$-invariant} if up to null sets we have $T_{g}(A) = A$ for every $g \in G$. The action of $G$ on $X$ is called \textit{ergodic} in case $\newmu(A) \in \{0,1\}$ for every $G$-invariant set $A \subseteq X$. In the following we shall write $\text{Fin}(G)$ for the set of finite subsets of $G$. Given a finite partition $\mathcal{P}$ of $X$ and a set $F \in \text{Fin}(G)$ we shall denote by \pagebreak
	\[H_{\newmu}^{F}(\mathcal{P}) := H_{\newmu}\left(\bigvee_{g \in F}T_{g}^{-1}(\mathcal{P})\right)\]
	the entropy of the refinement of $\mathcal{P}$ along $F$. \\
	
	\noindent\textit{Two commuting automorphisms.} We shall first discuss the case of two commuting automorphisms, which corresponds to an action of $\mathbb{Z}^{2}$. In this case a canonical symmetric set of generators is given by $\Gamma := \{\pm e_{i}\colon i=1,2\}$, where $e_{i}$ denotes the $i$-th standard basis vector of $\mathbb{Z}^{2}$. Let $\Xi$ denote the shift space $\Gamma^{\mathbb{N}}$ and consider the Bernoulli measure $\overline{\newnu}$ on $\Xi$ arising from the equidistribution on $\Gamma$. This corresponds to a simple recurrent random walk on $\mathbb{Z}^{2}$. It is easy to see that $\overline{\newnu}$ is $\newsigma$-invariant and ergodic. 
	
	Let $\{T_{g}\colon g \in \mathbb{Z}^{2}\}$ be an ergodic measure preserving action of $\mathbb{Z}^{2}$ on some standard probability space $(X,\newmu)$. Then, by Bufetov's criterion, the  corresponding step skew product $T$ is also ergodic, so we obtain an ergodic bundle RDS $(\Xi \times X,\overline{\newnu} \otimes \newmu,T)$. Now let $\mathcal{P}$ be a finite partition of $X$. Denoting by
	\[R_{n}(\newalpha) := \left\{\sum_{j=0}^{i-1}\newalpha_{j}\colon i= 0,...,n{-}1\right\}\]
	the range of the random walk up to time $n$ we may write 
	\begin{equation*}
		\begin{split}
			\frac{1}{n}H_{\newmu}^{n}(\mathcal{P})(\newalpha) &= \frac{|R_{n}(\newalpha)|}{n}\cdot\frac{1}{|R_{n}(\newalpha)|}~H_{\newmu}^{R_{n}(\newalpha)}(\mathcal{P}) \\ &\le \frac{|R_{n}(\newalpha)|}{n}~\frac{\log|\mathcal{P}|^{|R_{n}(\newalpha)|}}{|R_{n}(\newalpha)|} \\
			&= \frac{|R_{n}(\newalpha)|}{n}~\log|\mathcal{P}|.
		\end{split}
	\end{equation*}
	Now by a well known property of recurrent random walks, cf. \cite[Ch. I.4]{Spi76}, we obtain
	\[\lim_{n \to \infinity}\frac{\mathrel{\raisebox{-.4ex}{$1$}}}{\mathrel{\raisebox{+.1ex}{$n$}}}|R_{n}(\newalpha)| = 0\]
	and therefore
	\[\lim_{n \to \infinity}\frac{\mathrel{\raisebox{-.4ex}{$1$}}}{\mathrel{\raisebox{+.1ex}{$n$}}}H_{\newmu}^{n}(\mathcal{P})(\newalpha) = 0\]
	for $\newnu$-almost all $\newalpha \in \Xi$. Note that the left hand side is bounded by $\log|\mathcal{P}|$, so we may apply dominated convergence to obtain $h_{\newnu,\newmu}(\mathcal{P}) = 0$. Since $\mathcal{P}$ was arbitrary, this implies that the fiber entropy of the RDS is zero. \\
	
	\noindent\textit{Two non commuting automorphisms.} As a second example we consider the case of two non commuting automorphisms, which corresponds to an action of the free group $\mathbb{F}_{2}$ over the symmetric generator set $\Gamma := \{a^{\pm1},b^{\pm1}\}$. Denoting again by $\Xi$ the shift space $\Gamma^{\mathbb{N}}$ we consider the Markov measure $\newnu^{*}$ on $\Xi$ corresponding to the equidistributed initial vector $\pi$ on $\Gamma$ and the transition matrix $\Pi$ with entries given by
	\[\Pi(s,s') := \begin{cases} \frac{1}{3},~s' \neq s^{-1} \\ 0,~s' = s^{-1}\end{cases}\]
	for $s,s' \in \Gamma$. The construction guarantees that a set $[\textv]$ has positive probability if and only if the word $\textv$ is uncancellable over $\Gamma$. Thus, considering the product topology on $\Xi$, the measure $\newnu^{*}$ is supported on the subshift $\Xi_{0} \subseteq \Xi$ generated by the uncancellable words. It is not difficult to check that $\Pi$ is irreducible and $\pi$ satisfies the equation $\pi^{T}\Pi = \pi^{T}$. Accordingly $\newnu^{*}$ defines a $\newsigma$-invariant and ergodic measure on $\Xi$.
	
	For a given alphabet $\Lambda$ consider the space $\Lambda^{\mathbb{F}_{2}}$ equipped with the $\newsigma$-algebra generated by cylinder sets. Let $\overline{\newmu} = \overline{\pi}^{\mathbb{F}_{2}}$ be the Bernoulli measure on $\Lambda^{\mathbb{F}_{2}}$ corresponding to the equidistribution $\overline{\pi}$ on $\Lambda$ and consider the shift action $\{S_{g}\colon g \in \mathbb{F}_{2}\}$ of $\mathbb{F}_{2}$ on the probability space $\left(\Lambda^{\mathbb{F}_{2}},\overline{\newmu}\right)$ defined by
	\[S_{g}(\newbeta) := (\newbeta_{g^{-1}h})_{h \in \mathbb{F}_{2}}\]
	for $\newbeta\in \Lambda^{\mathbb{F}_{2}}$ and
	$g \in \mathbb{F}_{2}$. It is well known that this action is measure preserving and ergodic, see e.\@ g. \cite[Ch. 2.3.1]{KL16}. Furthermore, for all $s,s' \in \Gamma$ we find some $q \in \Gamma$ such that $\Pi(q,s) > 0$ and $\Pi(q,s') > 0$, so $\newnu^{*}$ satisfies the condition of Bufetov's criterion. Denoting by $S$ the corresponding step skew product, we may therefore conclude that the bundle RDS $(\Xi \otimes \Lambda^{\mathbb{F}_{2}}, \newnu^{*} \otimes \overline{\newmu},S)$ is ergodic.
	
	Let $\overline{\mathcal{N}}$ denote the partition of $\Lambda^{\mathbb{F}_{2}}$ consisting of cylinder sets of the form $\big\{\newbeta_{e} = \lambda\big\}$ for $\lambda \in \Lambda$. It is not difficult to see that
	\[H_{\overline{\newmu}}^{F}(\overline{\mathcal{N}}) = |F|H_{\overline{\newmu}}~\!\big(\overline{\mathcal{N}}\big) = |F|\log|\Lambda|\]
	for every $F \in \text{Fin}(\mathbb{F}_{2})$, so we have
	\[\inf_{F \in \text{Fin}\left(\mathbb{F}_{2}\right)}\frac{1}{|F|}H_{\overline{\newmu}}^{F}(\overline{\mathcal{N}}) = \log|\Lambda| > 0.\]
	By a theorem of Bowen, cf. \cite[Thm. 2.13]{Bow20}, and the non-amenability of $\mathbb{F}_{2}$ this implies
	\begin{align} \label{eq24}
		\sup_{\mathcal{P}}\inf_{F \in \text{Fin}\left(\mathbb{F}_{2}\right)}\frac{1}{|F|}H_{\overline{\newmu}}^{F}(\mathcal{P}) = \mathlarger{\mathlarger{\infinity}}.
	\end{align}
	For $\newalpha \in \Xi$ and $n \ge 1$ we consider the finite subsets $F_{n}(\newalpha)$ of $\mathbb{F}_{2}$ given by
	\[F_{n}(\newalpha) := \left\{\prod_{j=1}^{i}\newalpha_{i-j}\colon i = 0,...,n{-}1\right\}.\]
	Since $\newnu^{*}$ is supported on $\Xi_{0}$ and there is a one-to-one correspondence between the uncancellable words over $\Gamma$ and the elements of $\mathbb{F}_{2}$, the increasing initial segments of $\newnu$-almost all sequences $\newalpha$ correspond to pairwise different elements of $\mathbb{F}_{2}$. Considering the bijection on $\mathbb{F}_{2}$ mapping an element $s_{j}\cdot ... \cdot s_{0}$ to the mirrored element $s_{0}\cdot ... \cdot s_{j}$ we obtain therefore $|F_{n}(\newalpha)| = n$ for $\newnu$-almost all $\newalpha \in \Xi$. Thus we have
	\begin{align*}
		\frac{1}{n}H_{\overline{\newmu}}^{n}\left(\mathcal{P}\right)(\newalpha) &= \frac{1}{|F_{n}(\newalpha)|}H_{\overline{\newmu}}^{F_{n}(\newalpha)}(\mathcal{P}) \\ &\ge \inf_{F \in \text{Fin}\left(\mathbb{F}_{2}\right)}\frac{1}{|F|}H_{\overline{\newmu}}^{F}(\mathcal{P})
	\end{align*}
	for $\newnu$-almost all $\newalpha \in \Xi$, all $n \ge 1$ and any finite partition $\mathcal{P}$ of $X$. Integrating and taking the limit yields
	\[h_{\newnu,\newmu}(\mathcal{P}) \ge \inf_{F \in \text{Fin}\left(\mathbb{F}_{2}\right)}\frac{1}{|F|}H_{\overline{\newmu}}^{F}(\mathcal{P}).\]
	Thus, by (\ref{eq24}) the fiber entropy $h_{\newnu,\newmu}$ of the RDS is infinite. Note that in view of the random Kolmogorov-Sinai theorem discussed above and the fact that the fiber entropy with respect to a partition is finite this implies that there can be no random generator for the above system.
	
	\section{Proof of the main theorem} \label{Proofs}
	
	Let $(\Xi \times X, \newnu \otimes \newmu, T)$ be a finite bundle RDS and $\mathcal{P}$ be a finite partition of $X$. As indicated in the introduction we shall define the \textit{upper} and \textit{lower conditional orbit complexity} $\overline{\mathcal{C}}_{\mathcal{P}}(\textx|\newalpha)$ and $\underline{\mathcal{C}}_{\mathcal{P}}(\textx|\newalpha)$ of a state $\textx \in X$ along a sequence $\newalpha \in \Xi$ as the upper and lower conditional algorithmic complexity of the $\mathcal{P}$-$\newalpha$-name of $\textx$ relative to $\newalpha$. We have already noted that this definition is meaningful for $\newnu\otimes\newmu$-almost all pairs $(\newalpha,\textx) \in \Xi \times X$. Furthermore, using (\ref{Equ6}) it is not difficult to see that the above complexities are invariant under the skew product, i.\@ e. they depend only on the orbit of $(\newalpha,\textx)$ under $T$, which justifies the naming.
	
	A deterministic dynamical system $(X,\newmu,T)$ can be identified with the special case of a $1$-ary alphabet $\Theta = \{\newtheta\}$ with $\newnu$ being the trivial probability measure on $\Xi = \{\newtheta^{\infty}\}$. It is not difficult to see that in this case one has $h_{\newnu,\newmu}(\mathcal{P}) = h_{\newmu}(\mathcal{P})$. Furthermore, the $\mathcal{P}$-$\newtheta^{\infty}$-name of $\textx$ equals the $\mathcal{P}$-name of $\textx$ and, by the properties of conditional complexities discussed above, we obtain $\overline{\mathcal{C}}_{\mathcal{P}}(\textx|\newtheta^{\infty}) = \overline{\mathcal{C}}_{\mathcal{P}}(\textx)$ as well as $\underline{\mathcal{C}}_{\mathcal{P}}(\textx|\newtheta^{\infty}) = \underline{\mathcal{C}}_{\mathcal{P}}(\textx)$ for $\newmu$-almost all $\textx \in X$. Thus, Theorem \ref{main} contains Brudno's classical theorem as a special case.
	
	We will now turn to the proof of Theorem \ref{main}. As we shall see it suffices to restrict to the case of random shifts. To this end let $\Theta$ and $\Lambda$ be alphabets and let $\Omega$ denote the corresponding configuration space $\Lambda^{(\Theta^{*})}$. For $\lambda \in \Lambda$ let $C_{\lambda}$ be the cylinder set consisting of all $\newbeta \in \Omega$ such that $\newbeta_{\mathbf{e}} = \lambda$. We shall call the partition $\mathcal{N}^{*}$ of $\Omega$ consisting of the sets $C_{\lambda}$ the \textit{natural partition} of $\Omega$. For $\textu \in \Theta^{k}$ and $\textv \in \Lambda^{k}$ we define the measurable set $[\textu|\textv] \subseteq \Omega$ by
	\[[\textu|\textv] := \bigcap_{i=0}^{k-1} S_{u,i}^{-1}(C_{v_{i}}),\]
	where $S_{u,0}$ denotes the identity and $S_{u,i}$ denotes the map $S_{u_{i-1}} \circ ... \circ S_{u_{0}}$. (The specification of the last symbol of $\textu$ in the notation may seem superfluous. However, this notation will turn out to be very useful below.) In the following we shall often write $\newmu[\textu|\textv]$ for the value of a measure $\newmu$ on the set $[\textu|\textv]$. Obviously the following statement is a special case of Theorem \ref{main}\@.	
	\begin{theorem} \label{mainshift}
		Let $(\Xi \times \Omega, \newnu \otimes \newmu,S)$ be an ergodic random shift and $\mathcal{N}^{*}$ be the natural partition of $\Omega$. Then for $\newnu\otimes\newmu$-almost all $(\newalpha,\newbeta) \in \Xi \times \Omega$ we have
		\[\overline{\mathcal{C}}_{\mathcal{N}^{*}}(\newbeta|\newalpha) = \underline{\mathcal{C}}_{\mathcal{N}^{*}}(\newbeta|\newalpha) = {h}_{\newnu,\newmu}(\mathcal{N}^{*}).\]
	\end{theorem}
	
	Moreover, it will turn out that Theorem \ref{mainshift} is actually equivalent to Theorem \ref{main}. To see this consider an arbitrary finite ergodic bundle RDS $(\Xi \times X, \newnu \otimes \newmu,T)$ together with a finite partition $\mathcal{P}$ of $X$. Let $\Omega = \mathcal{P}^{(\Theta^{*})}$ be the configuration space corresponding to the alphabets $\Theta$ and $\mathcal{P}$. Denoting by $T_{u}$ the transformation $T_{u_{|u|-1}} \circ ... \circ T_{u_{0}}$ for $\textu \in \Theta^{*}$, where $T_{\mathbf{e}} = \text{Id}$, we may consider the map $\psi\colon X \to \Omega$ sending a state $\textx \in X$ to the a.\@ s. unique element $\newbeta \in \Omega$ satisfying $T_{u}(\textx) \in \newbeta_{u}$ for every $\textu \in \Theta^{*}$. It is easy to check that $\psi$ is measurable. Furthermore we obtain $\psi \circ T_{\newtheta} = S_{\newtheta} \circ \psi$ for all $\newtheta \in \Theta$. Denoting by $\psi_{*}\newmu$ the push-forward measure of $\newmu$ under $\psi$ this implies
	\begin{equation*}
		\begin{split}
			\psi_{*}\newmu\big(S_{\newtheta}^{-1}(B)\big) &=\newmu\big(\psi^{-1}S_{\newtheta}^{-1}(B)\big) = \newmu\big(T_{\newtheta}^{-1}\psi^{-1}(B)\big) \\ &= \newmu\big(\psi^{-1}(B)\big) = \psi_{*}\newmu(B)
		\end{split}
	\end{equation*}
	for all measurable sets $B \subseteq \Omega$. Thus $\psi_{*}\newmu$ is invariant under $S_{\newtheta}$ for all $\newtheta \in \Theta$. Accordingly we obtain a random shift $(\Xi \times \Omega,\newnu \otimes \psi_{*}\newmu,S)$.
	
	Next, consider the map $\Phi = \text{Id} \otimes \psi$ and note that $\Phi_{*}(\newnu\otimes\newmu) = \newnu \otimes \psi_{*}\newmu$. By the above observation $\Phi$ commutes with the skew-products $T$ and $S$, which implies that for an $S$-invariant set $C$ the preimage $\Phi^{-1}(C)$ is $T$-invariant, so we get
	\[\newnu \otimes \psi_{*}\newmu(C) = \Phi_{*}(\newnu\otimes\newmu)(C) = \newnu\otimes\newmu\big(\Phi^{-1}(C)\big) \in \{0,1\}\]
	by the ergodicity of $T$. Accordingly $S$ is an ergodic transformation. Denoting for $\textu \in \newtheta^{*}$ and $i \le |\textu|$ by $T_{u,i}$ the transformations $T_{u_{i-1}} \circ ... \circ T_{u_{0}}$, where $T_{u,0} = \text{Id}$, we obtain furthermore
	
	\begin{equation*}
		\begin{split}
			H_{\newnu,\psi_{*}\newmu}^{n}(\mathcal{N}^{*}) &= -\sum_{|u| = n} \newnu[\textu]\sum_{|v| = n}\psi_{*}\newmu\Big([\textu|\textv]\Big)~\log\psi_{*}\newmu\Big([\textu|\textv]\Big)
			\\ &= -\sum_{|u| = n} \newnu[\textu]\sum_{|v| = n}\newmu\Big(\psi^{-1}[\textu|\textv]\Big)~\log \newmu\Big(\psi^{-1}[\textu|\textv]\Big) \\
			&= -\sum_{|u|=n}\newnu[\textu] \sum_{|v|=n}\newmu\left(\bigcap_{i=0}^{n-1}T_{u,i}^{-1}\big(\textv_{i}\big)\right)~\log\newmu\left(\bigcap_{i=0}^{n-1}T_{u,i}^{-1}\big(\textv_{i}\big)\right) \\
			&= ~~H_{\newnu,\newmu}^{n}(\mathcal{P})
		\end{split}
	\end{equation*}
	for every $n \ge 1$, which gives $h_{\newnu,\newmu}(\mathcal{P}) = h_{\newnu,\psi_{*}\newmu}(\mathcal{N}^{*})$. Thus, by Theorem \ref{mainshift}, we find a set $D$ with $\newnu\otimes\psi_{*}\newmu(D) = 1$ such that $\overline{\mathcal{C}}_{\mathcal{N}^{*}}(\newbeta|\newalpha) = \underline{\mathcal{C}}_{\mathcal{N}^{*}}(\newbeta|\newalpha) = h_{\newnu,\newmu}(\mathcal{P})$ for all $(\newalpha,\newbeta) \in D$. Noting that the $\mathcal{N}^{*}$-$\newalpha$-name of $\psi(\textx)$ equals the $\mathcal{P}$-$\newalpha$-name of $\textx$ (modulo a relabelling of the symbols in the underlying alphabet) we obtain therefore
	\begin{equation} \label{Equ11}
		\overline{\mathcal{C}}_{\mathcal{P}}(\textx|\newalpha) = \overline{\mathcal{C}}_{\mathcal{N}^{*}}\big(\psi(\textx)|\newalpha\big) = h_{\newnu,\newmu}(\mathcal{P}) = \underline{\mathcal{C}}_{\mathcal{N}^{*}}\big(\psi(\textx)|\newalpha\big) = \underline{\mathcal{C}}_{\mathcal{P}}(\textx|\newalpha)
	\end{equation}
	for all $(\newalpha,\textx) \in \Phi^{-1}(D)$. Since by definition we have  	
	\[\newnu\otimes\newmu\big(\Phi^{-1}(D)\big) = \Phi_{*}(\newnu\otimes\newnu)(D) =  \newnu\otimes\psi_{*}\newmu(D) = 1,\]
	this shows that (\ref{Equ11}) is actually valid for $\newnu \otimes \newmu$-almost all $(\newalpha,\textx) \in \Xi \times X$. This verifies that Theorem \ref{mainshift} implies Theorem \ref{main}.
	
	Consequently we are reduced to show Theorem \ref{mainshift}\@ . We shall split the proof into two lemmas, which together will imply the assertion. To this end we have to introduce a further notion. Let $A \subseteq \Lambda^{*}$ be a finite set of words over an alphabet $\Lambda$. By a \textit{prefix free code} of $A$ we mean an injective map $\kappa\colon A \to \{0,1\}^{*}$ such that the set of codewords $\kappa(A)$ is prefix free. A useful tool for constructing prefix free codes is Kraft's inequality, which states that for a sequence $r_{1},...,r_{m}$ of natural numbers the existence of a prefix free set of words $\textv^{1},...,\textv^{m} \in \{0,1\}^{*}$ with lengths $|\textv^{i}| = r_{i}$ is equivalent to the condition that the powers $2^{-r_{1}},...,2^{-r_{m}}$ sum up to a number smaller or equal than $1$, see \cite[Thm. 2.1.2]{Rom92}.
	
	\begin{lemma}
		Let $(\Xi \times \Omega, \newnu \otimes \newmu,S)$ be an ergodic random shift and let $\mathcal{N}^{*}$ denote the natural partition of $\Omega$. Then for $\newnu\otimes\newmu$-almost all $(\newalpha,\newbeta) \in \Xi \times \Omega$ we have
		\[\overline{\mathcal{C}}_{\mathcal{N}^{*}}(\newbeta|\newalpha) \le h_{\newnu,\newmu}(\mathcal{N}^{*}).\]
	\end{lemma}
	
	\textit{Proof:} We will show that for any fixed $k \ge 1$ and $\newnu \otimes\newmu$-almost all $(\newalpha,\newbeta)$ we obtain
	\begin{equation} \label{Equ12}
		\overline{\mathcal{C}}_{\mathcal{N}^{*}}(\newbeta|\newalpha) ~\!\le~\! \frac{\mathrel{\raisebox{-.4ex}{$1$}}}{k}~\!\! H_{\newnu,\newmu}^{k}(\mathcal{N}^{*}) + \frac{\mathrel{\raisebox{-.4ex}{$1$}}}{k}.
	\end{equation}
	Since the right hand side converges to $h_{\newnu,\newmu}(\mathcal{N}^{*})$ as $k \to~\!\!\!\! \mathrel{\raisebox{-.2ex}{$\mathlarger{\mathlarger{\infty}}$}}$, this will suffice to prove the assertion.  \pagebreak
	
	To this end fix $\newbeta \in \Omega$ and $\newalpha \in \Xi$ and let $\newomega$ denote the $\mathcal{N}^{*}$-$\newalpha$-name of $\newbeta$. (Note that here \textit{every} $\newbeta$ has a unique $\mathcal{N}^{*}$-$\newalpha$-name.) In the following we will think of the sequences $\newomega$ and $\newalpha$ as decomposed into blocks of length $k$, where we shall denote by $\newomega^{i}$ and $\newalpha^{i}$ the $i$-th $k$-blocks $\newomega_{\smash{(i-1)k}}...\newomega_{\smash{ik-1}}$ and $\newalpha_{(i-1)k}...\newalpha_{ik-1}$ of $\newomega$ and $\newalpha$ respectively. We want to code the $k$-blocks of $\newomega$ in dependence of the corresponding $k$-blocks of $\newalpha$. Note that for every $\textu \in \Theta^{k}$ the probabilities $\newmu[\textu|\textv]$ with $\textv \in \Lambda^{k}$ sum up to $1$, so by Kraft's inequality we obtain  a prefix free code $\kappa_{u}$ of $\Lambda^{k}$ such that
	\begin{equation} \label{Equ13}
		|\kappa_{u}(\textv)| \le -\log\newmu[\textu|\textv] + 1
	\end{equation}
	for every $\textv \in \Lambda^{k}$. This allows us to code a $k$-block $\newomega^{i}$ by the code word $\kappa_{\newalpha^{i}}(\newomega^{i})$. Now consider an oracle machine $M$, which for every $m \ge 1$ restores the initial segment  $\newomega^{1}...\newomega^{m}$ of $\newomega$ from the concatenation of codewords
	\[\kappa_{\newalpha^{1}}(\newomega^{1})....\kappa_{\newalpha^{m}}(\newomega^{m}),\]
	when provided with $\newalpha$ as oracle. Such a machine could work as follows. Beginning with $i=1$ it scans the input symbol per symbol unless the scanned prefix matches a codeword of the form $\kappa_{\newalpha^{i}}(\textv)$, where it knows the block $\newalpha_{i}$ from the oracle. In this case it writes $\textv$ at the end of the output and iterates the procedure with the remaining part of the input and the updated counter $i{+}1$. Since the number of possible $k$-blocks in $\newomega$ and $\newalpha$ is finite, the possible codes can be stored in a finite table (encoded in the states of the machine in a suitable way). Furthermore, since the set of possible codewords is prefix free in every step, there can be no mismatch. To secure that the machine terminates only on those inputs, which are of the required form, we may make it go into an infinite loop, whenever it reaches the end of the input without having found an appropriate codeword.
	
	Next consider a window of length $k$ sliding simultaneously through $\newalpha$ and $\newomega$ with a stepwidth of $d$ and let $\mathbf{a}_{d,m}^{u,v}(\newalpha,\newbeta)$ denote the frequency of the pair $(\textu,\textv)$ under the first $m$ scans. Then, by definition, the value $m~ \mathbf{a}_{k,m}^{u,v}(\newalpha,\newbeta)$ gives the number of occurrences of $(\textu,\textv)$ under the first $m$ $k$-blocks, so we have
	\[\sum_{|u|,|v|=k} \!\!\! m~\!\mathbf{a}_{k,m}^{u,v}(\newalpha,\newbeta) = m.\]
	Consequently, by (\ref{Equ13}), we obtain
	\begin{equation}\label{Equ14}
		\begin{split}
			C_{M}(\newomega^{1}...\newomega^{m}|\newalpha) &\le \sum_{|u|,|v|=k} \!\!\!m~\!\mathbf{a}_{k,m}^{u,v}(\newalpha,\newbeta)~\!|\kappa_{u}(\textv)| \\
			&\le - \sum_{|u|,|v| = k} \!\!m~\mathbf{a}_{k,m}^{u,v}(\newalpha,\newbeta)\log\newmu[\textu|\textv] + m.
		\end{split}
	\end{equation}
	Setting
	\[\widehat{H}_{m}^{k}(\newalpha,\newbeta) := - \sum_{|u|,|v|=k} \mathbf{a}_{k,m}^{u,v}(\newalpha,\newbeta)\log\newmu[\textu|\textv]\]
	for $m \ge 1$ (\ref{Equ14}) together with (\ref{Equ3}) yields
	\begin{equation} \label{Equ15}
		\begin{split}
			C_{mk}(\newomega|\newalpha) &\le C_{M}(\newomega^{1}...\newomega^{m}|\newalpha) + O(1) \\ &\le m~\!\widehat{H}_{m}^{k}(\newalpha,\newbeta) + m + O(1).
		\end{split}
	\end{equation}
	It seems natural here to think of $\mathbf{a}_{k,m}^{u,v}(\newalpha,\newbeta)$ as an estimator of $\newnu[\textu]\newmu[\textu|\textv]$ and representing it in the form
	\[\mathbf{a}_{k,m}^{u,v}(\newalpha,\newbeta) = \frac{1}{m}\sum_{i=0}^{m-1}\mathbf{1}_{[u]\times[u|v]}\big(\smash{S}^{ki}(\newalpha,\newbeta)\big)\]
	we are tempted to write
	\[\lim_{m \to \infinity}\mathbf{a}_{k,m}^{u,v}(\newalpha,\newbeta) = \mathbb{E}\mathbf{1}_{[u]\times[u|v]} = \newnu[\textu]\newmu[\textu|\textv]\]
	and therefore $\widehat{H}_{m}^{k}(\newalpha,\newbeta) \to H_{\newnu,\newmu}^{k}(\mathcal{N}^{*})$ for $m \to \mathlarger{\mathlarger{\infinity}}$ and $\newnu\otimes \newmu$-almost all $(\newalpha,\newbeta) \in \Xi \times \Omega$ in the light of Birkhoff's ergodic theorem. Unfortunately we cannot justify this here directly. Although the power $\smash{S}^{k}$ of the skew-product $S$ is measure preserving and $\smash{\mathbf{a}_{k,m}^{u,v}(\newalpha,\newbeta)}$ as well as $\smash{\widehat{H}_{m}^{k}(\newalpha,\newbeta)}$ will thus indeed converge for $\newnu \otimes \newmu$-almost all $(\newalpha,\newbeta)$, $\smash{S}^{k}$ will in general not be ergodic, so we may not infer that the limit equals $\mathbb{E}\mathbf{1}_{[u]\times[u|v]}$ as we would like to do. (Note that under the stronger assumption that $S$ is weakly mixing we could draw this conclusion.)
	
	However, with a little trick we obtain the desired result in a way that suffices for our purposes. The crucial observation is that
	\begin{equation*}
		\begin{split}
			\mathbf{a}_{1,mk}^{u,v}(\newalpha,\newbeta) &= \frac{\mathrel{\raisebox{-.4ex}{$1$}}}{km}\sum_{i=0}^{km-1}\mathbf{1}_{[u]\times[u|v]}\left(\smash{S}^{i}(\newalpha,\newbeta)\right) \\
			&= \frac{\mathrel{\raisebox{-.4ex}{$1$}}}{km}\sum_{r=0}^{k-1}\sum_{j=0}^{m-1}\mathbf{1}_{[u]\times[u|v]}\left(\smash{S}^{r+jk}(\newalpha,\newbeta)\right) \\
			&=  \frac{\mathrel{\raisebox{-.4ex}{$1$}}}{k}\sum_{r=0}^{k-1}\mathbf{a}_{k,m}^{u,v}\big(\smash{S}^{r}(\newalpha,\newbeta)\big).
		\end{split}
	\end{equation*}
	So evaluating $\widehat{H}_{m}^{k}$ on the tuples $(\newalpha,\newbeta), S(\newalpha,\newbeta), ...,  \smash{S}^{k-1}(\newalpha,\newbeta)$ and averaging over the values yields
	\[\frac{1}{k}\sum_{r=0}^{k-1}\widehat{H}_{m}^{k}\big(\smash{S}^{r}(\newalpha,\newbeta)\big) = - \sum_{|u|,|v| = k} \mathbf{a}_{1,mk}^{u,v}(\newalpha,\newbeta)\log\newmu[\textu|\textv].\]
	Since by assumption $S$ is ergodic, we may now apply Birkhoff's theorem to obtain that
	\begin{align*}
		\lim_{m \to \infinity}\mathbf{a}_{1,mk}^{u,v}(\newalpha,\newbeta) &= \lim_{m \to \infinity}\frac{\mathrel{\raisebox{-.4ex}{$1$}}}{mk}\sum_{i=0}^{mk-1}\mathbf{1}_{[u]\times[u|v]}\big(\smash{S}^{i}(\newalpha,\newbeta)\big) \\
		&= \mathbb{E}\mathbf{1}_{[u]\times[u|v]} = \newnu[\textu]\newmu[\textu|\textv]
	\end{align*}
	for $\newnu\otimes\newmu$-almost all $(\newalpha,\newbeta) \in \Xi \times \Omega$. Since the limit \[\widehat{H}^{k}(\newalpha,\newbeta) := \lim_{n \to \infinity}\widehat{H}_{m}^{k}(\newalpha,\newbeta)\]
	exists $\newnu\otimes\newmu$-almost surely as mentioned above and $S$ is measure preserving, we may conclude that the limits $\widehat{H}^{k}(\newalpha,\newbeta)$,  $\widehat{H}^{k}\big(S(\newalpha,\newbeta)\big)$, ...  $\widehat{H}^{k}\big(S^{k-1}(\newalpha,\newbeta)\big)$ exist simultaneously and we have
	\begin{equation*}
		\begin{split}
			\frac{\mathrel{\raisebox{-.4ex}{$1$}}}{k}\sum_{r=0}^{k-1}\widehat{H}^{k}\big(\smash{S}^{r}(\newalpha,\newbeta)\big) &= \lim_{m \to \infinity}\frac{\mathrel{\raisebox{-.4ex}{$1$}}}{k}\sum_{r=0}^{k-1}\widehat{H}_{m}^{k}\big(\smash{S}^{r}(\newalpha,\newbeta)\big) \\
			&= -\sum_{|u|=k}\newnu[\textu]\sum_{|v|=k}\newmu[\textu|\textv]\log\newmu[\textu|\textv] \\
			&= ~H_{\newnu,\newmu}^{k}(\mathcal{N}^{*})
		\end{split}
	\end{equation*}
	for $\newnu\otimes\newmu$-almost all $(\newalpha,\newbeta) \in \Xi \times \Omega$. Thus for every such pair $(\newalpha,\newbeta)$ there is at least one $r \! \in \! \{0,...,k{-}1\}$ such that
	\begin{equation} \label{Equ16}
		\widehat{H}^{k}\big(\smash{S}^{r}(\newalpha,\newbeta)\big) \le H_{\newnu,\newmu}^{k}(\mathcal{N}^{*}).
	\end{equation}
	Fixing such a pair $(\newalpha,\newbeta)$ together with the respective $r$ let $(\newgamma,\neweta)$ denote the pair $S^{r}(\newalpha,\newbeta)$ and let $\newtau$ denote the $\mathcal{N}^{*}$-$\newgamma$-name of $\neweta$. For an arbitrary $n \ge 1$ write $n = mk + \ell$ with $m,k, \ell \in \mathbb{N}$ and $\ell < k$. Then, by (\ref{Equ5}) together with (\ref{Equ2}) and the fact that the conditional complexity is bounded from above by the plain complexity modulo some constant, we obtain
	
	\begin{equation*}
		\begin{split}
			C_{n}(\newtau|\newgamma) &\le C_{mk}(\newtau|\newgamma) + C(\newtau_{mk}...\newtau_{n-1}|\newgamma) + 2\log (mk^{2}) + O(1) \\
			&\le C_{mk}(\newtau|\newgamma) + \ell \log|\Lambda| + 4\log (mk) + O(1)
		\end{split}
	\end{equation*}
	and therefore
	\[\overline{\mathcal{C}}(\newtau|\newgamma) = \limsup_{n \to \infinity}\frac{\mathrel{\raisebox{-.4ex}{$1$}}}{n}C_{n}(\newtau|\newgamma)~\!\!\le~\!\!\limsup_{m \to \infinity}\frac{\mathrel{\raisebox{-.4ex}{$1$}}}{mk}~\!\!C_{mk}\big(\newtau|\newgamma\big),\]
	which by (\ref{Equ15}) and (\ref{Equ16}) implies  \[\overline{\mathcal{C}}_{\mathcal{N}^{*}}(\neweta|\newgamma)~\!\!\le~\!\! \frac{\mathrel{\raisebox{-.4ex}{$1$}}}{k}~\!\!\widehat{H}^{k}\big(\smash{S}^{r}(\newalpha,\newbeta)\big) + \frac{\mathrel{\raisebox{-.4ex}{$1$}}}{k}~\!\!\le~\!\! \frac{\mathrel{\raisebox{-.4ex}{$1$}}}{k}~\!\! H_{\newnu,\newmu}^{k}(\mathcal{N}^{*}) + \frac{\mathrel{\raisebox{-.4ex}{$1$}}}{k}.\]
	Finally, since the upper conditional orbit complexity only depends on the orbit of $(\newalpha,\newbeta)$ under $S$, the values of $\overline{\mathcal{C}}_{\mathcal{N}^{*}}(\neweta|\newgamma)$ and $\overline{\mathcal{C}}_{\mathcal{N}^{*}}(\newbeta|\newalpha)$ coincide, so the above argument shows the validity of (\ref{Equ12}) for $\newnu\otimes\newmu$-almost all $(\newalpha,\newbeta) \in \Xi \times \Omega$. \hfill $\Box$
	
	\begin{lemma}
		Let $(\Xi \times \Omega, \newnu \otimes \newmu,S)$ be an ergodic random shift and let $\mathcal{N}^{*}$ denote the natural partition of $\Omega$. Then for $\newnu\otimes\newmu$-almost all $(\newalpha,\newbeta) \in \Xi \times \Omega$ we have
		\[\underline{\mathcal{C}}_{\mathcal{N}^{*}}(\newbeta|\newalpha) \ge h_{\newnu,\newmu}(\mathcal{N}^{*}).\]
	\end{lemma}
	
	\textit{Proof:} For a given pair $(\newalpha,\newbeta) \in \Xi \times \Omega$ we denote by $\newbeta_{\newalpha,i}$ the inscription of $\newbeta$ at the coordinate $\newalpha_{0}...\newalpha_{i-1}$ and define the conditional complexities
	\[K_{n}(\newalpha,\newbeta) := K(\newbeta_{\newalpha,0}...\newbeta_{\newalpha,n-1}|\newalpha)\]
	for $n \ge 1$. After relabelling the symbols of the underlying alphabet we may identify the sequence $(\newbeta_{\newalpha,n})_{n \in \mathbb{N}}$ with the $\mathcal{N}^{*}$-$\newalpha$-name of $\newbeta$. We shall show that for $\newnu\otimes\newmu$-almost all $(\newalpha,\newbeta)$ there exists an $N \ge 1$ such that we have
	\[K_{n}(\newalpha,\newbeta) > J^{~\!\!\!n}_{\newmu}(\mathcal{N}^{*})(\newalpha,\newbeta) - 2\log n\]
	for all $n \ge N$. Since by the above observation the limes inferior of the left hand side divided by $n$ equals $\underline{\mathcal{C}}_{\mathcal{N}^{*}}(\newbeta|\newalpha)$ and the right hand side divided by $n$ converges to $h_{\newnu,\newmu}(\mathcal{N}^{*})$ almost surely by the random SMB theorem, this will suffice to prove the assertion.
	
	To this end let $E$ be the set of all $(\newalpha,\newbeta)$, where the above assertion fails. For a fixed $\newalpha \in \Xi$ let $E^{\newalpha}$ denote the $\newalpha$-fiber $\big\{\newbeta \in \Omega\colon(\newalpha,\newbeta) \in E\big\}$ of $E$. Defining for $n \ge 1$ the sets
	\[E_{n}^{\newalpha} := \Big\{\newbeta \in \Omega\colon K_{n}(\newalpha,\newbeta) \le J_{\newmu}^{~\!\!\! n}(\mathcal{N}^{*})(\newalpha,\newbeta) - 2\log n\Big\}\]
	we may represent $E^{\newalpha}$ in the form
	\[E^{\newalpha} = \bigcap_{N=1}^{\infinity}\bigcup_{n \ge N}E_{n}^{\newalpha}.\]
	Introducing for $\newbeta \in \Omega$ and $\textv \in \Lambda^{*}$ the sets
	\[Q_{\newalpha}^{v} := \bigcap_{i=0}^{|v|-1}S_{\newalpha,i}^{-1}\big(C_{v_{i}}\big)\]
	and
	\[Q^{n}_{\newalpha,\newbeta} := \bigcap_{i=0}^{|v|-1}S_{\newalpha,i}^{-1}\big(C_{\newbeta_{\newalpha,i}}\big)\]
	and noting that
	\[J_{\newmu}^{~\!\!\! n}(\mathcal{N}^{*})(\newalpha,\newbeta) = -\log\newmu\big(Q^{n}_{\newalpha,\newbeta}\big)\]
	we obtain
	\[\newbeta \in E_{n}^{\newalpha} \Leftrightarrow \newmu\big(Q^{n}_{\newalpha,\newbeta}\big) \le n^{-2}~2^{-K_{n}(\newalpha,\newbeta)}.\]
	Furthermore, writing $D$ for the set of all $\textv \in \Lambda^{n}$ such that $Q_{\newalpha}^{v} \cap E_{n}^{\newalpha} \neq \emptyset$, we have obviously \[E_{n}^{\newalpha} \subseteq \bigcup_{v \in D}Q^{v}_{\newalpha}.\]
	Note that for all $\textv \in D$ we find some $\newbeta \in E_{n}^{\newalpha}$ such that $\textv = \newbeta_{\newalpha,0}...\newbeta_{\newalpha,n-1}$ and therefore $Q^{v}_{\newalpha} = Q^{n}_{\newalpha,\newbeta}$, so we have $\newmu\big(Q^{v}_{\newalpha}\big) \le n^{-2}2^{-K(v|\newalpha)}$. Accordingly we obtain the upper bound
	\[\newmu\big(E_{n}^{\newalpha}\big) \le \sum_{v \in D}\newmu\big(Q_{\newalpha}^{v}\big) \le n^{-2}\sum_{v \in D}2^{-K(v|\newalpha)}.\]
	Since by definition the numbers $K(\textv|\newalpha)$, $\textv \in D$, are the lengths of a set of prefix free words, Kraft's inequality implies 	
	\[\sum_{v \in D}2^{-K(v|\newalpha)} \le 1.\]
	Using the fact that $\sum_{n=1}^{\infinity}n^{-2} < \mathlarger{\mathlarger{\infinity}}$ an application of the Borel-Cantelli lemma yields that $E^{\newalpha}$ is a null set. Since $\newalpha$ was arbitrary, this implies
	\[\newnu \otimes \newmu(E) = \int\newmu\big(E^{\newalpha}\big)~d\newnu(\newalpha) = 0\]
	by Cavalieri's principle, which proves the claim. \hfill $\Box$ \\
	
	It remains to show Theorem \ref{AR-comp}, which is now easily obtained. Consider a finite bundle RDS $(\Xi \times X,\newnu \otimes \newmu,T)$ and let $\mathcal{P}$ be a finite partition of $X$. Then $\mathcal{N} \times \mathcal{P}$ is a finite partition of $\Xi \times X$. Note that $\overline{\mathcal{C}}_{\mathcal{N}}(\newalpha) = \overline{\mathcal{C}}(\newalpha)$ and $\underline{\mathcal{C}}_{\mathcal{N}}(\newalpha) = \underline{\mathcal{C}}(\newalpha)$ and recall that we have $h_{\newnu} = h_{\newnu}(\mathcal{N})$. Thus, since $\newsigma$ and $T$ are ergodic, the values $\mathcal{C}(\newalpha)$, $\mathcal{C}_{\mathcal{N}\times\mathcal{P}}(\newalpha,\textx)$ and $\mathcal{C}_{\mathcal{P}}(\textx|\newalpha)$ exist by Theorem \ref{Brudno} and Theorem \ref{main} and satisfy
	\[\mathcal{C}_{\mathcal{P}}(\textx|\newalpha) = h_{\newnu,\newmu}(\mathcal{P}) = h_{\newnu\otimes\newmu}(\mathcal{N} \times \mathcal{P}) - h_{\newnu} = \mathcal{C}_{\mathcal{N}\times\mathcal{P}}(\newalpha,\textx) - \mathcal{C}(\newalpha)\]
	for $\newnu \otimes \newmu$-almost all $(\newalpha,\textx) \in \Xi \times X$ by (\ref{Equ9}).
	
	As a final remark we note that under the above assumptions the (unconditional) algorithmic complexity of the $\mathcal{P}$-$\newalpha$-name of a state $\textx \in X$ along a sequence $\newalpha \in \Xi$ can also be quantified for $\newnu \otimes \newmu$-almost all $(\newalpha,\textx) \in \Xi \times X$. It is not difficult to see that the upper and lower algorithmic complexity of the $\mathcal{P}$-$\newalpha$-name of $\textx$ equals the upper and lower algorithmic complexity of the $\{\Xi\} \times \mathcal{P}$-name of $(\newalpha,\textx)$ with respect to $T$. Therefore, as a consequence of Brudno's classical theorem, both values coincide $\newnu \otimes \newmu$-almost surely and equal the entropy $h_{\newnu \otimes \newmu}(\{\Xi\} \times \mathcal{P})$ of $T$ with respect to the partition $\{\Xi\} \times \mathcal{P}$.
	\section*{Acknowledgments} The results of this paper extend results of the author's diplo\-ma thesis. The author thanks his advisor Felix Pogorzelski for freely sharing his ideas on the topic and giving valuable hints in many stages of the work as well as for constant motivation and substantial support during the writing process. Moreover he thanks the anonymous reviewers for their useful comments and suggestions on an earlier version of the manuscript. Finally he gratefully acknowledges financial support through a grant of the German-Israeli Foundation for Scientific Research and Development (GIF).

\end{document}